\documentclass[12pt]{amsart}
\usepackage{amsmath}
\usepackage{enumerate}
\usepackage{amssymb}
\usepackage{amsbsy}
\usepackage{amsfonts}
\usepackage[arrow]{xy}
\usepackage{diagrams}
\theoremstyle{plain}
\newtheorem{thm}{Theorem}[section]
\newtheorem{prop}[thm]{Proposition}
\newtheorem{lem}[thm]{Lemma}

\theoremstyle{definition}
\newtheorem{dfn}[thm]{Definition}

\newtheorem{rmk}[thm]{Remark}
\numberwithin{equation}{section}
\newcommand{\sm}{\left(\begin{smallmatrix}}
\newcommand{\esm}{\end{smallmatrix}\right)}

\newfont{\FieldFont}{msbm10 scaled\magstep1}

\newcommand{\R}{{\mathcal{R}}}

\newcommand{\Z}{{\mathbb{Z}} }

\newcommand{\pf}{\noindent\bf Proof }

\def\MAT#1#2#3#4{\Bigl(
\begin{matrix}#1&#2\\#3&#4\end{matrix}\Bigr)}
\def\SMA#1#2#3#4{\bigl(
\begin{smallmatrix}#1&#2\\#3&#4\end{smallmatrix}\bigr)}

\newcommand{\beq}{\begin{eqnarray*}}
\newcommand{\eeq}{\end{eqnarray*}}

 \def\MOD#1{\;(\text{mod}\;#1)}

\newcommand{\leb}{\left[}
\newcommand{\rib}{\right]}

\begin{document}

\title{  Periods of
Jacobi forms  and Hecke operator }

\author{YoungJu  Choie }

 \address{Department of Mathematics and PMI\\
 Pohang University of Science and Technology\\
 Pohang, 790--784, Korea}
 \email{yjc@postech.ac.kr}

 \author{Seokho Jin }

  \address{Department of Mathematics and PMI\\
  Pohang University of Science and Technology\\
  Pohang, 790--784, Korea}
  \email{archimed@postech.ac.kr }

 \thanks{Keynote: Jacobi form, Hecke Operator, period}
 \thanks{This work
 was partially supported by NRF 2012047640, NRF 2011-0008928 and  
NRF 2008-0061325
}

\begin{abstract}
A Hecke action on  the space of periods of   cusp forms,
 which is compatible with that  on the space of   cusp forms, 
was first computed using continued fraction\cite{M0} and 
  an explicit algebraic formula of  Hecke operators
acting on the space of period functions of   modular forms was derived
by studying   the rational period functions\cite{CZ}. As an application
an elementary proof of the Eichler-Selberg trace formula
 was derived\cite{Z2}.
 Similar modification has been applied to period space of 
Maass cusp forms with spectral parameter $s$\cite{Mu, Mu1, M01}. 
In this paper we study the space of period functions of Jacobi forms by 
means of Jacobi integral and  give an explicit
description of Hecke operator  acting on this space.
a Jacobi Eisenstein series  $E_{2,1}(\tau,z)$ of weight $2$ and index $1$ 
is discussed  as an example.
 Periods of Jacobi integrals are already appeared as a disguised form 
 in the work of Zwegers to study   Mordell integral coming 
 from Lerch sums\cite{Ze} and mock Jacobi forms are 
 typical example of Jacobi integral\cite{Z3}.
 
\end{abstract}
 \maketitle

%\today

%%%%%%%%%%%%%%%%%%%%%%%%%%%%%%%%%%%%%%%%%%%%%
%%%%%%%%%%%%%%%%%%%%%%%%%%%%%%%%%%%%%%%%%%%%%%%

%%%%%%%%%%%%%%%%%%%%%%%%%%%%%%%%%%%%%%%%%%%%%%%%%%%%%%%%%%%%%%%%
\section{\bf{Introduction}}

%%%%%%%%%%%%%%%%%%%%%%%%%%%%%%%%%%%%%%%%%%%%%%%%%%%%%%%%%%%%%

Period functions of   modular forms
 have played  an important role 
 to understand the arithmetics on cusp forms \cite{KZ}.
Manin\cite{M0} studied  a Hecke action on the space of period  of   cusp forms, 
which is compatible with that  on the space of   cusp forms, in terms of continued fractions.
    Later  an explicit algebraic description of a Hecke operator  
on the space of period functions of 
 modular forms was given by studying  
 the rational period functions\cite{CZ}.
  As an application,  
a new elementary proof of the Eichler-Selberg trace formula was derived \cite{Z2}.   
Moreover,  
various modifications of period theory   have been also developed.
A notion of rational period functions
 has been introduced and
 completely classified
in the case of cofinite subgroups of $SL_{2}(\mathbb{Z})$
 \cite{K1, K2, K0, A}.
Some period functions  were
attached bijectively to Maass cusp forms 
according to the spectral parameter $s$   and
its  cohomological counterpart was described (see \cite{LZ, BLZ}).
 Similar modification to Hecke operators
  on the period space has also been applied to that of 
Maass cusp forms with spectral parameter $s$ \cite{Mu, Mu1, M01}. 

Historically, Eichler\cite{E} and Shimura\cite{S0} discovered an isomorphism
between a space of cusp forms and   Eichler cohomology group,
attaching period polynomials to cusp forms.
A notion of Eichler integral for
arbitrary real weight with multiplier system has been introduced  
and shown that there is an isomorphism between the space of modular forms
and   Eichler cohomology group\cite{K1}.
Along this line, recently  
a notion of Jacobi integral has been introduced \cite{CL}
and    shown   there is also an isomorphism
between the space of Jacobi cusp forms and
the corresponding Eichler cohomology group\cite{ChL}.
 Mock Jacobi form\cite{Z3} is one of the typical example of Jacobi integral. 

The main purpose of this paper is to   give an explicit algebraic description of 
Hecke operator  on  the period functions
attached to the Jacobi integrals.
This is an   
analogous result of \cite{CZ} to the case of Jacobi forms. 

This paper is organized as follows: in section $2$ we state the main results and
 section $3$ discusses about a
Jacobi 
Eisenstein series $E_{2,1}(\tau,z)$  as an example of Jacobi integral with period functions.  
Basic definitions and notations are given in 
section $4$ and section $5$ gives   proofs of the main theorems 
by introducing various   properties which are modified 
from the results about rational period functions  \cite{CZ}.
%%%%%%%%%%%%%%%%%%%%%%%%%%%%%%%%%%%%%%%%%%

\section{\bf{Statement of Main Results}}

%%%%%%%%%%%%%%%%%%%%%%%%%%%%%%%%%%%%%%%%%

Let $f$ be an element of $J_{k,m}^{\int}(\Gamma(1)),$  that is, $f$ is 
a real analytic function
$f : \mathbb{H}\times \mathbb{C} \rightarrow \mathbb{C}$
satisfying a certain growth
condition with the following functional equation,
\begin{equation}\label{p}
(f|_{k,m}\gamma)(\tau,z)
=f(\tau,z)+ P_{\gamma}(\tau,z) ,
\forall \gamma\in \Gamma(1)^J,
\end{equation}
with   $P_{\gamma}\in \mathcal{P}_{ m},$
where $\mathcal{P}_{ m}$ is a set of holomorphic functions   
\begin{equation}
 \mathcal{P}_{ m}=\{g : \mathcal{H}\times\mathbb{C}
 \rightarrow \mathbb{C} | \,
 |g(\tau,z)| < K(|\tau|^\rho+v^{-\sigma})
 e^{2\pi m\frac{
y^2}{v} },\  \mbox{ for some   $K,\rho, \sigma >0$}\} 
\end{equation}
(  $v=Im(\tau)$ and $y=Im(z)$ ).\\
$P_{\gamma}$ in (\ref{p}) is called
  a {\bf {period function}} of $f .$
If $P_{\gamma}(\tau,z)=0,
\forall \gamma \in \Gamma(1)^J,$     $f$ is a usual Jacobi form (see \cite{EZ}).
It turns out that each element of the following set 
\begin{equation}\label{per}
\mathcal{P}er_{k,m}:= \{ P:\mathcal{H}\times
\mathbb{C} \rightarrow \mathbb{C} \,  | \,
\sum_{j=0}^{3} P|_{k,m}T^j
=\sum_{j=0}^{5}P|_{k,m}U^j=0 \}  \
\end{equation}
( $T=[\sm 0 & -1 \\ 1 & 0 \esm, (1, 0)], S=[\sm 1 & 1 \\ 0 & 1 \esm, (0,0)],$
$U=ST $) \\
generates a system of period functions $\{P_{\gamma}|\gamma \in \Gamma(1)^J\}$ of  $f\in J_{k,m}^{\int}(\Gamma(1)).$

We also consider the following   
subspace:

\begin{equation}\label{ei}
EJ^{\int}_{k,m}(\Gamma(1)) :=
\{f \in J_{k,m}^{\int}(\Gamma(1))  | \,
f|_{k,m}[I, (1,0)]=f \}. \end{equation}

It turns out that the following set  
\begin{equation}\label{per}
\mathcal{EP}er_{k,m}:= \{ P:\mathcal{H}\times
\mathbb{C} \rightarrow \mathbb{C} \,  | \,
\sum_{j=0}^{3} P|_{k,m}T^j
=\sum_{j=0}^{5}P|_{k,m}U^j=
P - P|_{k,m}[-I,(1,0)]=0 \}.
\end{equation}
is a generating set of a system  of
period functions $\{P_{\gamma} | \, \gamma\in \Gamma(1)^J\}$ of $f \in  EJ_{k,m}^{\int}(\Gamma(1)).$

For each positive integer $n,$ define two Hecke operators $\mathcal{V}_{n}^{\infty}$ 
and  $ \mathcal{T}_{n}^{\infty}$ by

\begin{equation}\label{he1}
\left(f|_{k,m}\mathcal{V}_{n}^{\infty}\right)(\tau,z):
 =n^{k-1}\sum_{ad=n, a>o\atop b\MOD d}d^{-k}
 f\left(\frac{a\tau+b}{d},az\right)
\end{equation}

\begin{equation}\label{he2}
 (f|_{k,m}\mathcal{T}_{n}^{\infty})(\tau,z):=n^{k-4}
\sum_{\begin{array}{cc} ad=n^2,
\gcd(a,b,d)=\square \atop a,d>0, b\MOD d
,  \small{X,Y \in \mathbb{Z}/n\mathbb{Z}}\end{array}}
(f|_{k,m} [\sm a & b \\ 0 & d \esm,
(X,Y)])(\tau,z),
\end{equation}

Then 
\begin{thm}\label{Hecke1}
\begin{enumerate}
\item  
$\left(f|_{k,m}\mathcal{V}_{n}^{\infty}\right)
 \in J_{k,mn}^{\int}(\Gamma(1))$ if   $f\in J_{k,m}^{\int}(\Gamma(1)).$ 

\item Let $f|_{k,m}T=f+P_T$ and
 $(f|_{k,m}\mathcal{V}_{n}^{\infty})|_{k,mn}T=
(f|_{k,m}\mathcal{V}_{n}^{\infty})
+ \hat{P_T} .$
Then
$$\hat{P_T}= n^{\frac{k}{2}-1}\biggl({P}_{T}|_{k,m}\biggl[\MAT
{\frac{1}{\sqrt{n}}}00{\frac{1}{\sqrt{n}}},
(0,0)\biggr]\biggr) |_{k,m}
\tilde{\mathcal{V}}_{n },
$$
where
\begin{eqnarray*}
\tilde{\mathcal{V}}_{n }&=&\sum_{ad-bc=n\atop a>c>0,
d>-b>0}\biggl\{\biggl[\MAT abcd,(0,0)\biggr]
+\biggl[\MAT a{-b}{-c}d,(0,0)\biggr]\biggr\}\\
&&+\sum_{ad=n\atop -\frac{1}{2}d<b\leq \frac{1}{2}d}
\biggl[\MAT ab0d,(0,0)\biggr]+
\sum_{ad=n\atop-\frac{1}{2}a
<c\leq\frac{1}{2}a, c\neq 0}\biggl[\MAT a0cd,(0,0)\biggr].
\end{eqnarray*}

\end{enumerate}
\end{thm}

\begin{thm}\label{Hecke2}

\begin{enumerate}
\item
If  $f \in EJ_{k,m}^{\int}(\Gamma(1)), $
then $f|_{k,m}\mathcal{T}_{n}^{\infty}
\in EJ_{k,m}^{\int}(\Gamma(1)).$

\item Let $f|_{k,m}T=f +P_{T}$ and
   $(f|_{k,m}\mathcal{T}_{n}^{\infty})
    |_{k,m}T=(f|_{k,m}\mathcal{T}_{n}^{\infty})+ \tilde{P_T}.$
Then  $$ \tilde{P_T}=n^{k-4} P_T|_{k,m}\tilde{\mathcal{T}_n},$$ 
where
\begin{eqnarray*}
 &&\tilde{\mathcal{T}}_n=\sum_{\tiny{
\begin{array}{ccc}
ad-bc=n^{2}, gcd(a,b,c,d)=\square\\ a>c>0,
 d>-b>0\\  X,Y\in\Z/n\Z\end{array}}}
\left\{[\sm a&b\\c&d\esm,(X,Y)]+[\sm a&
-b\\-c&d\esm, (X,Y)]\right\}\\
&&+\sum_{\tiny{\begin{array}{ccc}
ad=n^{2}, gcd(a,b,d)=\square\\ -\frac{1}{2}d<b\leq
\frac{1}{2}d\\ X,Y\in\Z/n\Z\end{array}}} [\sm a&b\\0&d\esm,(X,Y)]
+\sum_{\tiny{\begin{array}{ccc}ad=n^{2},gcd(a,c,d)
=\square\\ -\frac{1}{2}a<c\leq\frac{1}{2}a, c\neq0\\ X,Y\in \Z/n\Z\end{array}}}
[\sm
a& 0\\c&d\esm,(X,Y)].
\end{eqnarray*}

\end{enumerate}
\end{thm}

\begin{rmk}
\begin{enumerate}
\item  There are   Hecke operators    acting
on the space of Jacobi forms \cite{EZ}.
  One needs to choose a special set of representatives
to apply
Hecke operators
to Jacobi integral $f$
since $f$ is not $\Gamma(1)^J$-invariant.

\item  Note that $\mathcal{T}_n^{\infty}$ acts only on the subspace
 $ EJ_{k,m}^{\int}(\Gamma(1)) $
(see section \ref{Hec} for details).

\item
It is shown that there is an isomorphism
 between the space of Jacobi cusp forms and  Eichler cohomology
 group with some coefficient module(there, denoted by 
 $\mathcal{P}_{\mathcal{M}}^{e}$) corresponding to 
 period functions for $f\in EJ_{k,m}^{\int}(\Gamma(1)).$  
 Hence the Hecke operator $\tilde{\mathcal{T}_{n}}$ can 
 also be regarded as an Hecke operator on the Eichler cohomology group.
 \end{enumerate}
\end{rmk}

\medskip

%%%%%%%%%%%%%%%%%%%%%%%%%%%%%%%%%%%%%%%%%%
\section{\bf{ Example of Jacobi Integral  }}

%%%%%%%%%%%%%%%%%%%%%%%%%%%%%%%%%%%%%%%%%%
%%%%%%%%%%%%%%%%%%%%%%%%%%%%%%%%%%%%%%%%%%

It is well known that Eisenstein series
$E_2(\tau):=1-24\sum_{n\geq 1}( \sum_{0<d|n} d) q^n$
is not a modular form on $\Gamma(1)$, but it is a modular integral(see \cite{K2}).
Similarly Jacobi-Eisenstien series $E_{2,1}(\tau,z)$
 is not a Jacobi form  on $\Gamma(1)$,
but it is a Jacobi integral as we explain below.

The following  Jacobi-Eisenstein series was studied in \cite{C2}:
$$E_{2,1}(\tau,z) =
-12\sum_{n,r\in\mathbb{Z},
 r^{2}\leq 4n}H(4n-r^{2})q^{n}\zeta^{r},
 q=e^{2\pi i \tau}, \zeta=e^{2 \pi i z}, $$
where $H(n)$ be the class number of the quadratic forms of
discriminant $-n$ with $H(0)=-\frac{1}{12}.$

The theta decomposition of $ E_{2,1}(\tau,z)$ is

$$E_{2,1}(\tau,z)
=-12\cdot{\binom{\mathcal{H}_{0}(\tau)}
{\mathcal{H}_{1}(\tau)}}^{t}\cdot
\binom{\vartheta_{1,0}(\tau,z)}{\vartheta_{1,1}(\tau,z)},  $$
where
 $$\mathcal{H}_{\mu}(\tau):=\sum_{N\geq 0, N\equiv -\mu^{2}\MOD 4 }
 H(N )q^{\frac{N}{4} },\
 \vartheta_{1,\mu}(\tau,z)=\sum_{r\in\mathbb{Z}\atop r\equiv
\mu\MOD{2 }}q^{\frac{r^{2}}{4}}\zeta^{r}.$$

\begin{rmk}
\begin{enumerate}\item $E_{2,1}(\tau,z)$ was defined in
\cite{C2} as the holomorphic part of
$$E_{2,1}^{*}(\tau,z;s)=\frac{1}{2}
\sum_{c,d\in\mathbb{Z},gcd(c,d)=1}
\sum_{\lambda\in\mathbb{Z}}
\frac{e^{2\pi i(\lambda^{2}\frac{a\tau+b}{c\tau+d}
+2\lambda\frac{z}{c\tau+d}-\frac{cz^{2}}{c\tau+d})}}
{(c\tau+d)^{2}|c\tau+d|^{2s}}$$
at $s=0$.

\item A correspondence among $E_2(\tau), E_{2,1}(\tau,z)$ 
and $\mathcal{H}_{\frac{3}{2}}(\tau):=\sum_{n\geq 0}H(n)q^n$ 
was studied in \cite{C2}.\\
\item
$ \mathcal{H}_{\frac{3}{2}}(\tau)=\mathcal{H}_0(4\tau)+\mathcal{H}_1(4\tau)$
since $H(N)=0$ unless $N\equiv 0,3\MOD 4$.
\end{enumerate}

\end{rmk}

The following function $\mathcal{F}$
 transforms like a modular form of weight
$\frac{3}{2}$ on the group $\Gamma_{0}(4) $ (see \cite{HZ}, p 91-92):
  $$\mathcal{F}(\tau)=\mathcal{H}_{\frac{3}{2}}(\tau)
+v^{-\frac{1}{2}}\sum_{\ell=
-\infty}^{\infty}
\beta(4\pi l^{2}v)q^{-\ell^{2}},$$
where
$
 \beta(x):=\frac{1}{16\pi}
 \int_{1}^{\infty}u^{-\frac{3}{2}}e^{-xu}du ( x\geqq 0 )  $ 
and 
\begin{eqnarray*}
 \binom{v^{-\frac{1}{2}}\sum_{l \equiv 0 \MOD{2}}
\beta(\pi l^{2}v)q^{-\frac{l^{2}}{4}}}{v^{-\frac{1}{2}}\sum_{l \equiv 1 \MOD{2}}
\beta(\pi l^{2}v)q^{-\frac{l^{2}}{4}}}=\frac{1+i}{16\pi}\binom{\int_{-\bar{\tau}}^{i\infty}
(t+\tau)^{-\frac{3}{2}}\vartheta_{1,0}(t,0)dt}{\int_{-\bar{\tau}}^{i\infty}
(t+\tau)^{-\frac{3}{2}}\vartheta_{1,1}(t,0)dt}.
\end{eqnarray*}

Consider a function $$\varphi(\tau,z):=\mathcal{F}_{0}(\tau)\vartheta_{1,0}(\tau,z)
 +\mathcal{F}_{1}(\tau)\vartheta_{1,1}(\tau,z)
 =\binom{\mathcal{F}_{0}(\tau)}{\mathcal{F}_{1}(\tau)}^t
 \binom{\vartheta_{1,0}(\tau,z)}{\vartheta_{1,1}(\tau,z)},$$
 with
\begin{eqnarray*}
\mathcal{F}_{\mu}(\tau)=\mathcal{H}_{\mu}(\tau)
+v^{-\frac{1}{2}}\sum_{\ell \equiv \mu \MOD{2}}
\beta(\pi \ell^{2}v)q^{-\frac{\ell^{2}}{4}}, \mu=0,1.
\end{eqnarray*}
It is easy to check
  that $\varphi(\tau,z)$ transforms
  like a Jacobi form for $\Gamma(1)^{J}$ of weight $2$ and index $1$ and so
$
\varphi\big|_{2,1}
T=\varphi.
$ On the other hand, following the computation in \cite{HZ}, page 92 we see
$$\biggl\{\frac{1+i}{16\pi}\binom{\int_{-\bar{\tau}}^{i\infty}
(t+\tau)^{-\frac{3}{2}}\vartheta_{1,0}(t,0)dt}{\int_{-\bar{\tau}}^{i\infty}
(t+\tau)^{-\frac{3}{2}}\vartheta_{1,1}(t,0)dt}^t \cdot
\binom{\vartheta_{1,0}}{\vartheta_{1,1}}
\bigg|_{2,1}(T-E) \biggr\} (\tau,z)$$
$$= \frac{1+i}{16}\binom{\int_{i\infty}^{0}
(\tau+w)^{-\frac{3}{2}}\vartheta_{1,0}(w,0)dw}{\int_{i\infty}^{0}
(\tau+w)^{-\frac{3}{2}}\vartheta_{1,1}(w,0)dw}^{t}
\binom{\vartheta_{1,0}(\tau,z)}{\vartheta_{1,1}(\tau,z)},$$
using the transformation formula (for example, see \cite{C1}) 
of theta series $\theta_{1, \mu}(\tau,z).$
 So we conclude that 
$$(E_{2,1}|_{2,1}T)(\tau,z)=E_{2,1}(\tau,z)+P_{  E_{2,1}, T}(\tau,z),$$
where
$$P_{  E_{2,1}, T}(\tau,z)=  \frac{1+i}{16}\binom{\int^{i\infty}_{0}
(\tau+w)^{-\frac{3}{2}}\vartheta_{1,0}(w,0)dw}{\int^{i\infty}_{0}
 (\tau+w)^{-\frac{3}{2}}\vartheta_{1,1}(w,0)dw}^{t} \cdot
  \binom{\vartheta_{1,0}(\tau,z)}
 {\vartheta_{1,1}(\tau,z)}.
  $$

Next, for each prime $p,$  let
$$(E_{2}|_{2}\mathcal{T}_{p,2}^{\infty})(\tau) := p\sum_{ad=p, d>0\atop b\MOD d}
d^{-4}E_{2}(\frac{a\tau+b}{d}),$$
$$
(\mathcal{H}_{\frac{3}{2}}|_{\frac{3}{2}}\mathcal{T}_{p,
 \frac{3}{2}}^{\infty})(\tau) := \sum_{N\equiv 0,3\MOD 4}
\left( H(Np^{2})+\left( \frac{-N}{p}\right) H(N)+p
H\left(\frac{N}{p^{2}}\right)\right) q^{N},$$
and
$$
(E_{2,1}|_{2,1}\mathcal{T}_{p}^{\infty})(\tau,z) := 
p^{-2}\sum_{ad=p^{2},a,d>0\atop b\MOD d,
\gcd{(a,b,d)} = \square}
\sum_{\lambda,\mu\in \mathbb{Z}/p\mathbb{Z}}
(E_{2,1}|_{2,1}[\SMA a b 0 d,
(\lambda,\mu)])(\tau,z).
$$

It can be directly checked  that  the  following diagram commutes:

\begin{diagram}
E_{2}(\tau) & \rTo& (E_{2}|_{2}\mathcal{T}_{p,2}^{\infty})(\tau) \\
\uTo^{\varphi}  & \curvearrowright &  \uTo^{\varphi} \\
\mathcal{H}_{\frac{3}{2}}(\tau) & \rTo&  
(\mathcal{H}_{\frac{3}{2}}|_{\frac{3}{2}}\mathcal{T}_{p, \frac{3}{2}}^{\infty})(\tau)\\
\dTo^{\psi}  & \curvearrowright &  \dTo^{\psi} \\
E_{2,1}(\tau,z) & \rTo & (E_{2,1}|_{2,1}\mathcal{T}_{p}^{\infty})(\tau,z), \\
\end{diagram}
where $\varphi(\sum_{n\geq 0}c(n)q^{n}):=1-\frac{24}
{L(0,\left(\frac{D}{\cdot}\right))}\sum_{n\geq 1}
\sum_{d|n}\left(\frac{D}{d}\right)c\left(\frac{n^{2}}{d^{2}}|D|\right)q^{n}$
with  a fundamental
  discriminant $D$ 
and 
$\psi(\sum_{n\geq 0}c(n)q^{n}):=-12\sum_{n,r\in\mathbb{Z},
  r^{2}\leq 4n}c(4n-r^{2})q^{n}\zeta^{r} $ (see \cite{C2}, Theorem 3.2).

\begin{rmk}
There is a one-to-one correspondence,
which is Hecke equivariant,
  among
  the space of modular forms of weight $2k-2$($k$ even)
on  $\Gamma(1)$,
  the Kohnen plus space
  of weight $k-\frac{1}{2}$ on $\Gamma_{0}(4) $  and
the space of Jacobi
   forms of weight $k$ and index $1$ on
   $\Gamma(1)^{J}.$
The above diagram  shows that this
correspondence is extended to  the more general case
   $E_2(\tau), H_{\frac{3}{2}}(\tau)$ and $E_{2.1}(\tau,z).$
  \end{rmk}

The above commutative diagram
implies that 
$E_{2,1}|_{2,1}\mathcal{T}_{p}^{\infty} 
 =(p+1)E_{2,1} $ since
$E_{2}|_2 \mathcal{T}_{p,2}^{\infty} =(p+1)E_2 $
 (also see Knopp\cite{K2}). 
In summary  we have shown the following:

\begin{prop}\label{ES}

\begin{enumerate}
\item
 $(E_{2,1}|_{2,1}T)(\tau,z)=E_{2,1}(\tau,z)
+P_{E_{2,1},T}(\tau,z),$   where
$$P_{  E_{2,1}, T}(\tau,z)=  \frac{1+i}{16}\binom{\int^{i\infty}_{0}
(\tau+w)^{-\frac{3}{2}}\vartheta_{1,0}(w,0)dw}{\int^{i\infty}_{0}
 (\tau+w)^{-\frac{3}{2}}\vartheta_{1,1}(w,0)dw}^{t} \cdot
  \binom{\vartheta_{1,0}(\tau,z)}
 {\vartheta_{1,1}(\tau,z)}.
  $$

\item For each prime $p$,
let $(E_{2,1}|_{2,1} \mathcal{T}_p^{\infty}|_{2,1}T)(\tau,z)
=E_{2,1}(\tau,z)+\tilde{P}_{E_{2.1},T}(\tau,z).$
Then
$$\tilde{P}_{E_{2.1},T}(\tau,z)=p^{-2}P_{E_{2,1},T}|_{2,1}\tilde{\mathcal{T}}_p =
(p+1)\cdot P_{E_{2,1},T}(\tau,z).$$

\end{enumerate}
\end{prop}

%%%%%%%%%%%%%%%%%%%%%%%%%%%%%%%%%%%%%%%%%%%%%%%%%%%%%%%%%%%%
\section{\bf{Definitions and Notations}\label{Def}}
%%%%%%%%%%%%%%%%%%%%%%%%%%%%%%%%%%%%%%%%%%%%%%%%%%%%%%%%%

  Let  $\mathcal{H}$ be the
usual complex upper half plane, $\tau \in \mathcal{H}, z\in
\mathbb{C}$ and $ \tau=u+iv, z=x+iy, u,v,x,y\in \mathbb{R}.$
Take $k, m\in \mathbb{Z}.$  Let 
$$ \Gamma(1)^{J}:= \Gamma(1)  \ltimes \mathbb{Z}^{2 }= \{ [M,(\lambda,\mu)]| M
\in \Gamma(1), \lambda,\mu \in \mathbb{Z}  \}, (\Gamma(1)=SL_{2}(\mathbb{Z}))$$ be the full Jacobi group 
with  a group law
$$[M_1,(\lambda_1,\mu_1)][M_2,(\lambda_2,\mu_2)]=[M_1M_2, (\lambda ,\mu  )M_2
+(\lambda_2,\mu_2)]. $$ 

Let us introduce the following elements in $\Gamma(1)^J$: $ S =[\sm 1 & 1 \\ 0 &
1\esm ,(0,0)],
  T =[\sm 0 & -1 \\ 1 & 0 \esm,(1,0)],
T_0=[\sm 0 & -1 \\ 1 & 0 \esm,(0,0)],
U=ST=[  \sm 1 & -1 \\ 1 &
0\esm,(1,0)],$ 
$   I_{2}=[I,(1,0)], I_{1}=[I,(0,1)], E=[I, (0,0)].$

In \cite{C-J} it is known
 that $\Gamma(1)^J$ is generated by $S $ and $T $.
Also $\Gamma(1)^J$ is generated by $T $ and $U $ and 
 they satisfy the relations
$$T^4=U^6=E,$$
$$UT^{2}=T^{2}U[I,(-1,0)]=T^{2}[I,(0,-1)]U=[I,(0,1)]T^{2}U $$
and these are the defining relations for $\Gamma(1)^J.$

Also, let $G^{J}$ be the set of triples 
$[M,X,\zeta]$ $(M\in SL_{2}(\mathbb{R}),
 X\in\mathbb{R}^{2}, \zeta\in\mathbb{C}, |\zeta|=1)$. Then $G^{J}$ is a group via
$$[M,X,\zeta][M',X',\zeta']=[MM', XM'+X',
 \zeta\zeta'\cdot e^{2\pi i\det\binom{XM'}{X'}}],$$

acting on $\mathcal{H}\times
\mathbb{C}$ as 
$$\gamma (\tau,z)=\left(\frac{a\tau+b}{c\tau+d},
\frac{z+\lambda\tau+\mu}{c\tau+d}\right),   \gamma=[\sm a&b\\c&d\esm, (\lambda, \mu
),\zeta],$$ which defines a  usual slash operator on a function $f:
\mathcal{H}\times \mathbb{C}  \rightarrow \mathbb{C}$ defined by :
 $$(f|_{   k,m }\gamma) (\tau,z) :=j_{k,m}(\gamma,(\tau,z))
 f(\gamma(\tau,z)), $$
with
 $j_{k,m }(\gamma,(\tau,z)) := \zeta^{m}(c\tau+d)^{-k}e^{2 \pi im
 (- \frac{cz^2}{ c\tau+d}
 +\lambda^2 \tau+2\lambda
z+  \lambda\mu)}.$ 
Further let  $$f|_{k,m}[\SMA a b c d, (X,Y)]:
=f|_{k,m}[\SMA {a/\sqrt{\ell}}{b/\sqrt{\ell}}{c/\sqrt{l}}{d/\sqrt{l}},(X,Y),1],
 \mbox{if $ ad-bc=\ell>0$ }.$$ We will omit $\zeta$ if $\zeta=1$.

\begin{dfn}\label{integral}(Jacobi Integral)
\begin{enumerate}
\item A real analytic  periodic (in both variables)
 function $f: \mathcal{H} \times \mathbb{C}
\rightarrow \mathbb{C}$ is called a   {\bf{Jacobi Integral }} of
weight $  k\in  \mathbb{Z}$ and index $m \in \mathbb{Z}$ with
  {\bf{ period functions}}
$P_{\gamma}$ on $\Gamma(1)^{J}$ if it satisfies the following
relations:
\begin{enumerate}
\item[(i)] For all $\gamma \in \Gamma(1)^{J}$
\begin{equation}{\label{Int}}
f|_{  k, m }\gamma = f(\tau, z)+ P_{\gamma}(\tau,z),
 P_{\gamma}\in \mathcal{P}_{ m}.
\end{equation}
\item[(ii)] It satisfies a growth condition, when $v, y
\rightarrow \infty,$
$$|f(\tau,z)|v^{-\frac{k}{2}}
e^{-2\pi m \frac{y^2}{v} } \rightarrow 0.$$
\end{enumerate}
\end{enumerate}
 \end{dfn}
 The space of Jacobi integrals forms a vector space over
$\mathbb{C}$ and we denote it by $J_{  k, m}^{\int}(\Gamma(1)^{J}).$

The periodicity condition on $f$ is equivalent to saying that
 $f|_{ k,m}S=f|_{  k,m} I_1=f $  so that
 $P_{S}(\tau,z)= P_{I_1}(\tau,z)=0.$ 
 A set of  period functions
\begin{equation}\label{pset}
 \{P_{\gamma} | \gamma \in\Gamma(1)^J\} \mbox{\, of $f \in J_{k,m}^{\int}(\Gamma(1))$}
 \end{equation}
 satisfies the following consistency condition:
$$P_{\gamma_1 \gamma_2} = P_{\gamma_1 }|_{ k,m} \gamma_2 +
P_{\gamma_2},\ \text{for all}\ \gamma_1, \gamma_2\in \Gamma(1)^J.$$
So using the relations $T^4=U^6=E,$ it is easy to see that
$P_T$ satisfies
$$\sum_{j=0}^3P|_{k,m}T^j=\sum_{j=0}^5P|_{k,m}U^j=0$$ so that
$P_T$ belongs $\mathcal{P}er_{k,m}.$  
In fact it is shown in \cite{CL} that  $\mathcal{P}er_{k,m} $
 generates the space of period functions $\{P_{\gamma} | 
 \gamma\in \Gamma(1)^J\}$ of $f\in J_{k,m}^{\int}(\Gamma(1)).$  

Next we consider a subspace 
$$EJ_{k,m}^{\int}(\Gamma(1))=\{ g \in J_{k,m}^{\int}(\Gamma(1)) | \, 
g|_{k,m} I_2=g\}$$
By the similar method as that for $\mathcal{P}er_{k,m}$ 
we can show that  a set of period functions
 $\{P_{\gamma}|\gamma\in \Gamma(1)^J\}$ of $g$  is spanned by 
 $\mathcal{EP}er_{k,m}.$

\medskip

%%%%%%%%%%%%%%%%%%%%%%%%%%%%%%%%%%%%%%%%%%%%%%%%%%%%%%%%%%%%%%%%%%
\section{{\bf{Hecke Operators on the space of Jacobi Integral}}\label{Hec}}
%%%%%%%%%%%%%%%%%%%%%%%%%%%%%%%%%%%%%%%%%%%%%%%%%%%%%%%%%%%%%%%%%%
We start to prove Theorem \ref{Hecke2} and only
 sketch out the proof of Theorem \ref{Hecke1}
  briefly, since it is similar but simpler.

\subsection{Proof of Theorem \ref{Hecke2}}

It is immediate to see  $f|_{k,m}\mathcal{T}_{n}^{\infty}
\in EJ_{k,m}^{\int}(\Gamma(1)) $ if $f\in EJ_{k,m}^{\int}
(\Gamma(1))$ by checking the following:
$$
f|_{k,m}[\SMA {a/n}{(b+d)/n}0{d/n},(X,Y)] = f|_{k,m}[\SMA 1101, (0,0)]\cdot
[\SMA {a/n}{b/n}0{d/n}, (X,Y)],$$ 
$$f|_{k,m}[\SMA  {a/n}{b/n}0{d/n},(X+n,Y)] = f|_{k,m}[I,(0,-b)]\cdot[I,(d,0)]\cdot[\SMA  {a/n}{b/n}0{d/n}, (X,Y)],$$
$$f|_{k,m}[\SMA  {a/n}{b/n}0{d/n}, (X,Y+n)]  = f|_{k,m}[I,(0,a)]\cdot[\SMA  {a/n}{b/n}0{d/n}, (X,Y)].$$
 
 \medskip

To prove (2) in Theorem \ref{Hecke2}, first we  need to prove a couple of propositions and lemmata.  
For each integer
$n>0$ let ${\mathcal{M}}_n^J:
=\{\gamma=[\gamma_0, (X,Y)]  \,
| \, \gamma_{0}\in \frac{1}{n}(M_{2}(\Z)
/\{\pm 1\}), X,Y\in (\frac{\mathbb{Z}}{n})/n\mathbb{Z},
det(\gamma_{0})=1\}.$ Write ${\mathcal{M}}^J_{+}:
=\cup {\mathcal{M}}_n^J$  and $\mathcal{R}_n^J:=\Z[{\mathcal{M}}_n^J]$ and
$\mathcal{R}_{+}^J:=\Z[\mathcal{M}_{+}^J]=\oplus_n {\mathcal{M}}_n^J,$ for the
sets of finite integral linear combinations of elements of ${\mathcal{M}}_n^J$
and $\mathcal{M}_{+}^J, $ respectively. Then $\mathcal{R}_{+}^J$ is a
(non-commutative) ring with unity and is ``multiplicatively graded"
in the sense that $\mathcal{R}_{n}^J\mathcal{R}_{m}^J\subset \mathcal{R}_{mn}^J$
for all $m,n>0$; in particular, each $\mathcal{R}_{n}^J$ is a left
and right module over the group ring $\mathcal{R}_{1}^J=\Z[\Gamma(1)^J]$
of $\Gamma(1)^J$. Denote by $\mathcal{J}$ the right ideal
$$(1+T+T^{2}+T^{3})\mathcal{R}_{1}^J
+(1+U+U^{2}+U^{3}+U^{4}+U^{5})\mathcal{R}_{1}^J$$
of $\mathcal{R}_{1}^J$.
Finally denote $\mathcal{M}_{n,1}^{J}$ by
$\{\gamma\in\mathcal{M}_{n}^{J}|\ X,Y\in\mathbb{Z}/n\mathbb{Z}\}$. Then
\begin{prop}\label{Main} Let $\hat{\mathcal{T}}_{n}^{\infty}$ be $n^{-k+4}\mathcal{T}_n^{\infty}$. 
\begin{enumerate}\item 
For each integer $n\geq 1,$ $\hat{\mathcal{T}}_{n}^{\infty}(S-I)\equiv 0,\,
\hat{\mathcal{T}}_{n}^{\infty}(I_{1}-I)\equiv 0,
 \hat{\mathcal{T}}_{n}^{\infty}(I_{2}-I)
\equiv 0\MOD {(S-I)\R_n^J +(I_{1}-I)\R_n^J+(I_{2}-I)\R_{n}^{J}}$ and

$
\hat{\mathcal{T}}_{n}^{\infty}(T-I)\equiv
(T-I)\tilde{\mathcal{T}_n}  \pmod{(S-I)\R_n^J +(I_{1}-I)\R_n^J+(I_{2}-I)\R_{n}^{J}}
$
for a certain
element $\tilde{\mathcal{T}_n} \in \R_n^J.$
\item
This element is unique modulo $\mathcal{J}\R_n^J.$
\item
If $n$ and $n'$ are positive integers coprime to
$m$ then $\tilde{\mathcal{T}}_n,  \tilde{\mathcal{T}}_{n'}  \in \R_{+}^{J}$ satisfy
the product formula
$$\tilde{\mathcal{T}}_{n} \cdot \tilde{\mathcal{T}}_{n'}
=\sum_{d | \gcd{(n,n')}}d^{2k-3}
\tilde{\mathcal{T}}_{\frac{nn'}{d^2} }\pmod{(S-I)\R_{nn'}^J+
(I_{1}-I)\R_{nn'}^J+(I_{2}-I)\R_{nn'}^J}.$$
\end{enumerate}
\end{prop}

The proof of Proposition \ref{Main} is similar to that in
\cite{CZ} where the explicit Hecke operators 
on the space  rational period functions associated to
elliptic modular forms were computed.

{\pf of Proposition \ref{Main}} 
$(1)$ For the assertion
$\hat{\mathcal{T}}_{n}^{\infty}(S-I)\equiv 0$ we compute as in \cite{CZ}.
To claim $\hat{\mathcal{T}}_{n}^{\infty}(I_{1}-I)\equiv 0$ check that
$$\sum_{\begin{array}{cc} ad=n^2, \gcd(a,b,d)=\square\atop a,d>0, b\MOD d
, X,Y\in \mathbb{Z}/n\mathbb{Z}\end{array}}[\SMA{a/n}{b/n}{0}{d/n},(X,Y)][\SMA 1001, (0,1)]$$
$$= \sum_{\begin{array}{cc} ad=n^2, \gcd(a,b,d)=\square\atop a,d>0, b\MOD d
, X,Y\in \mathbb{Z}/n\mathbb{Z}\end{array}}[\SMA {a/n}{b/n}{0}{d/n},(X,Y+1)] $$
$$\equiv  \sum_{\begin{array}{cc} ad=n^2, \gcd(a,b,d)=\square\atop a,d>0, b\MOD d
, X,Y\in \mathbb{Z}/n\mathbb{Z}\end{array}}[\SMA {a/n}{b/n}{0}{d/n},(X,Y)].
$$
For the assertion $\hat{\mathcal{T}}_{n}^{\infty}(I_{2}-1)\equiv 0$ check that
$$
 \sum_{\begin{array}{cc} ad=n^2, 
 \gcd(a,b,d)=\square\atop a,d>0, b\MOD d
, X,Y\in \mathbb{Z}/n\mathbb{Z}\end{array}}[\SMA{a/n}{b/n}{0}{d/n},(X,Y)][\SMA 1001, (1,0)]$$
$$= \sum_{\begin{array}{cc} ad=n^2, \gcd(a,b,d)=\square\atop a,d>0, b\MOD d
, X,Y\in \mathbb{Z}/n\mathbb{Z}\end{array}}[\SMA{a/n}{b/n}{0}{d/n},(X+1,Y)]$$
$$\equiv  \sum_{\begin{array}{cc} ad=n^2, \gcd(a,b,d)=\square\atop a,d>0, b\MOD d
, X,Y\in \mathbb{Z}/n\mathbb{Z}\end{array}}[\SMA {a/n}{b/n}{0}{d/n},(X,Y)]. 
$$

To claim  $\hat{\mathcal{T}}_{n}^{\infty}(T-I)\equiv
(T-I)\tilde{\mathcal{T}_n}  \pmod{(S-I)\R_n^J +(I_{1}-I)\R_n^J+(I_{2}-I)\R_{n}^{J}}
$ we need the following lemma.

\begin{lem}\label{lemgen}
Take any $\gamma \in \Gamma(1)^J.$ Then
$$\gamma-I\in(T-I)\R_1^J+(S-I)\R_1^J.$$
\end{lem}
\noindent{\bf{Proof of Lemma \ref{lemgen}}}:
 This lemma follows by induction on the word length. Assume
this holds for some $\gamma\in \Gamma(1)^J.$ Note that
$T\gamma-I=(T-I)\gamma+(\gamma-I), $ $ S\gamma-I=(S-I)\gamma+(\gamma-I),
S^{-1}\gamma-I=(S-I)(-S^{-1}\gamma)+(\gamma-I)$ also belong to
$(T-I)\R_1^J+(S-I)\R_1^J$.  Since $S, T$ generates
$\Gamma(1)^J,$ lemma follows by induction on the word length.
{\qed}
\\

Now write $\hat{\mathcal{T}}_{n}^{\infty}$ as $\sum M_{i}$. Note that 
$$[\SMA {a/n}{b/n}{0}{d/n},(X,Y)]T=[\SMA{b/(b,d)}
{\beta}{d/(b,d)}{\delta},(0,0)][\SMA{(b,d)/n}
{-a\delta/n}{0}{n/(b,d)},(Y+1,-X)]$$ for some 
$\beta,\delta\in\mathbb{Z}$ such that $\frac{b}{(b,d)}\delta-\frac{d}{(b,d)}\beta=1$.
Hence for each index $i$
we can choose index $i'$ such that $M_{i}T
=\gamma_{i}M_{i'}$ for some $\gamma_{i}\in
 \Gamma(1)^{J}$ so that the set of $i'$'s 
 is equal to  that of $i$. Then $\hat{\mathcal{T}}_{n}^{\infty}(T-1)
 =\sum \gamma_{i}M_{i'}-M_{i}=\sum_{i}(\gamma_{i}-1)M_{i'}$, and this belongs to $(T-1)
\mathcal{R}_{n}^{J}+(S-1)\mathcal{R}_{n}^{J}$ by Lemma \ref{lemgen}.\\

$(2)$ To characterize the elements of $\mathcal{J}{\mathcal{M}}_n^J$   
consider an``\textit{acyclicity}" condition, which was introduced
in \cite{CZ} in the case of the modular group: suppose $V$ is an abelian group on which $\Gamma(1)^J$
acts on the left.  Then $V$ is a left $\R_1^J$-module.  For $X\in
\R_1^J$ let $Ker(X):=\{v\in V \, | \, Xv=0\}, Im(X):=\{Xv\, | \,
v\in V\}.$  We call $V$ {\em{acyclic}} if
$$Ker (1+T+T^{2}+T^{3}) \cap Ker(I+U+\cdots+U^5)=\{0\},$$
$$Ker(I-T)=Im(I+T+T^2+T^3),$$ and
$$
Ker(I-U)=Im(I+U+U^2+U^3+U^4+U^5).$$

Then the following holds:
\begin{lem}\label{acy1}
$\R_n^J$ is an acyclic $\R_1^J$-module for all $n$.
\end{lem}
{\pf of Lemma \ref{acy1}} First we claim that $Ker(1+T+T^{2}+T^{3}) \cap Ker(I+U+\cdots+U^5)=\{0\}:$ let $X=\sum
n_{\gamma}\gamma (n_{\gamma} \in \Z, \gamma\in {\mathcal{M}}_n^J)$ be an element of
$\R_n^J.$  Suppose that $X\in Ker(1+T+T^{2}+T^{3}) \cap Ker(I+U+\cdots+U^5).$ Take any
$r(\tau)=\frac{1}{\tau-a},$ where $a\in \mathbb{C}$ is not rational or
quadratic and let $q(\tau,z):=(r|_{k,m}X)(\tau,z).$ Then $q(\tau,z)$
behaves somewhat like the rational period functions in \cite{CZ},
 i.e., it has finite number of singularities as a function 
 of $\tau$(when $z$ is fixed), and these singularities are rational or real quadratic.
 But it can be seen that for some $z_{0}$, $q(\tau,z_{0})$ has a singularity at some point
$\tau=M^{-1}a$ with $n_{M}\neq 0$, contradicting to the fact that
$q(\tau,z_{0})$ have singularities only at rational
or quadratic irrational points.
On the other hand, if $X$ is left invariant under $T$,
then $n_{M}=n_{TM}=n_{T^{2}M}=n_{T^{3}M}$ for all $M$,
and since $M, TM, T^{2}M, $ and $T^{3}M$ are distinct,
this means that $X$ can be written as an integral linear combination
of elements $M+TM+T^{2}M+T^{3}M=(1+T+T^{2}+T^{3}) M\in\mathcal{R}_{N}^{J}$.
Similarly, $X=UX$ implies $n_{M}=n_{UM}=\cdots=n_{U^{5}M}$ for all $M$
and hence $X\in (1+U+\cdots+U^{5})\mathcal{R}_{n}^{J}$.
This proves the second hypothesis in the definition of acyclicity.

{\qed}

\begin{lem} \label{acy2} If $V$ is an acyclic $\Gamma-$module and $v\in V,$
then
$$(I-T)v\in (I-S)V  \Leftrightarrow v\in
(1+T+T^{2}+T^{3}) V+(I+U+U^2+U^3+U^4+U^5) V=\mathcal{J}V.$$
\end{lem}

{\pf of Lemma \ref{acy2}} The direction ``$\Leftarrow$" is true for any $\Gamma(1)^{J}$-module,
since $v=(1+T+T^{2}+T^{3}) x+(I+U+U^2+U^3+U^4+U^5) y$ implies
$$(I-T)v=(1-S^{-1}U)(I+U+\cdots+U^5) y=(S-1)S^{-1}(I+U+\cdots+U^5) y.$$
Conversely, assume that $V$ is acyclic and $(1-T)v=(1-S)w$ for some $w\in V$. Then
$$(1-T)(v-(1+T+T^{2})w)=(1-T)v-(1-T^{3})w=(T^{3}-S)w=(1-U)T^{3}w.$$
This element must vanish since $Im(1-T)\cap Im(1-U)\subset Ker (1+T+T^{2}+T^{3})
\cap Ker (I+U+\cdots+U^5)$=\{0\}. But then $v-(1+T+T^{2})w\in Im((1+T+T^{2}+T^{3}))$
and $T^{3}w\in Im((I+U+\cdots+U^5))$ by the second hypothesis
in the definition of acyclicity, so $v=(v-(1+T+T^{2})w)+(1+T+T^{2}+T^{3})
w-T^{3}w\in (1+T+T^{2}+T^{3}) V+(I+U+\cdots+U^5) V$.

{\qed}

Lemmas \ref{acy1} and \ref{acy2} give a characterization of
$\mathcal{J}\mathcal{R}_{n}^{J}$ as
$\{X\in \mathcal{R}_{n}^{J}|\ (1-T)X \in (1-S)\mathcal{R}_{n}^{J}\}$.

\medskip

$(2)$ The uniqueness of $\tilde{\mathcal{T}}_{n}$ modulo $\mathcal{J}\mathcal{R}_{n}^{J}$
follows immediately from this characterization
and the definition of $\tilde{\mathcal{T}}_{n}$. \\

$(3)$ Finally,
\begin{eqnarray*}
&&(T-1)(\tilde{\mathcal{T}}_n \cdot \tilde{\mathcal{T}}_{n'}
-\sum_{d | \gcd{(n,n')}}d^{2k-3}
\tilde{\mathcal{T}}_{\frac{nn'}{d^2} })\\
&\equiv& \hat{\mathcal{T}}_{n}^{\infty}(T-1)\tilde{\mathcal{T}}_{n'}-
\sum_{d | \gcd{(n,n')}}d^{2k-3}(T-1)\tilde{\mathcal{T}}_{\frac{nn'}{d^2} }\\
&\equiv& \hat{\mathcal{T}}_{n}^{\infty}[\hat{\mathcal{T}}_{n'}^{\infty}(T-1)-(S-1)X_{n'}
-(I_{1}-1)Y_{n'}-(I_{2}-1)Z_{n'}]-\sum_{d | \gcd{(n,n')}}
d^{2k-3}\hat{\mathcal{T}}_{\frac{nn'}{d^{2}}}^{\infty}(T-1)\\
&\equiv&(\hat{\mathcal{T}}_{n}^{\infty} \cdot \hat{\mathcal{T}}_{n'}^{\infty}
-\sum_{d | \gcd{(n,n')}}d^{2k-3}
\hat{\mathcal{T}}_{\frac{nn'}{d^{2}}}^{\infty})(T-1) \MOD{(S-1)
\mathcal{R}_{n}^{J}+(I_{1}-1)
\mathcal{R}_{n}^{J}+(I_{2}-1)\mathcal{R}_{n}^{J}},
\end{eqnarray*}
and
$$\hat{\mathcal{T}}_{n}^{\infty} \cdot \hat{\mathcal{T}}_{n'}^{\infty}
-\sum_{d | \gcd{(n,n')}}d^{2k-3}
\hat{\mathcal{T}}_{\frac{nn'}{d^{2}}}^{\infty}\equiv 0\MOD
{(S-1)\mathcal{R}_{n}^{J}
+(I_{1}-1)\mathcal{R}_{n}^{J}+(I_{2}-1)
\mathcal{R}_{n}^{J}}$$
by the usual calculation
for the commutation properties of
Hecke operators. This completes the proof of Proposition \ref{Main}. 

 Now we are ready to prove Theorem \ref{Hecke2}-(2):
\subsection{\pf of Theorem \ref{Hecke2}-(2)} For any $(X,Y)\in \Z^2/n\Z^2, $ the maps
$$A=[\sm a& b\\c& d \esm, (X,Y)] \rightarrow
I_{1}^{-1}TS^{-[\frac{a}{c}]}A=\left[\sm c & d\\-a+c[\frac{a}{c}] &
-b+d[\frac{a}{c}]\esm,(X,Y)\right]$$

$$B=[\sm a & b\\c& d \esm, (X,Y)] \rightarrow
I_{1}^{-1}S^{[\frac{d}{b}]}TB=\left[\sm -c+a[\frac{d}{b}]&
-d+b[\frac{d}{b}] \\ a& b\esm, (X,Y)\right],$$ where $[\,]$ denotes the
integral part, are to give inverse bijections between the sets
$$\mathcal{A}_n:=\{[\sm a& b\\ c& d\esm,(X,Y)] \in{\mathcal{M}}_{n,1}^J| a>c>0,
d>-b\geq0, b=0 \Rightarrow a\geq 2c, \gcd(a,b,c,d)=\square\}$$

and $$\mathcal{B}_n=\{[\sm a& b\\c & d\esm, (X,Y)]\in {\mathcal{M}}_{n,1}^J \, | \,
a>-c\geq 0, d>b>0, c=0 \Rightarrow d\geq 2b, \gcd(a,b,c,d)=\square\}.$$

Note that
$$\sum_{\tiny{\begin{array}{cc}ad-bc=n^{2},
gcd(a,b,c,d)=\square\\a>c>0,d>-b>0\\X,Y\in \mathbb{Z}/n\mathbb{Z}\end{array}}}\left\{
[\sm a&b\\c&d\esm,(X,Y)]-T[\sm a&-b\\-c&d\esm,(X,Y)]\right\}
$$
$$\equiv \sum_{\tiny{\begin{array}{cc}ad=n^{2}, gcd(a,b,d)=\square\\\frac{1}{2}d\geq b
>0\\X,Y\in \mathbb{Z}/n\mathbb{Z}\end{array}}}[\sm 0 & -d\\a& b\esm,(X,Y)]-
\sum_{\tiny{\begin{array}{cc}ad=n^{2}, gcd(a,c,d)=\square\\\frac{1}{2}a\geq c
>0\\X,Y\in \mathbb{Z}/n\mathbb{Z}\end{array}}}[\sm a & 0 \\c & d \esm,(X,Y)].$$

Conjugating this equation by $\alpha:=[\SMA {-1}00{1},(0,0)]$
changes the sign of all the off-diagonal coefficients
of each $2\times 2$ matrices and each $X$'s, and
preserves the property ``$\equiv$", since
$\alpha S\alpha^{-1}=[\SMA 1{-1}01,(0,0)],
\alpha I_{1}\alpha^{-1}=I_{1},$ and $\alpha I_{2}\alpha^{-1}=[\SMA 1001, (-1,0)]$.
Also since $X$ can be chosen freely in
$\mathbb{Z}/n\mathbb{Z}$, note that we may change $-X$ to $X$.
Adding the result to the original equation, we get
$$(I-T)\sum_{\tiny{
\begin{array}{cc}ad-bc=n^{2}, gcd(a,b,c,d)=\square\\a>c>0, d>-b>0\\ X,Y\in \mathbb{Z}/n\mathbb{Z}\end{array}}}
[(\sm a&b\\c&d\esm+\sm a&-b\\-c&d\esm),(X,Y)]$$
$$\equiv \sum_{\tiny{\begin{array}{cc}ad=n^{2}, gcd(a,b,d)=\square\\ 0<|b|\leq
\frac{1}{2}d\\X,Y\in \mathbb{Z}/n\mathbb{Z}\end{array}}} [(\sm 0 & -d\\ a& b
\esm-\sm d& 0
\\ b& a\esm ),(X,Y)].$$

Hence,

$$(I-T) \tilde{\mathcal{T}}_n\equiv
\sum_{\tiny{\begin{array}{cc}  ad=n^{2},
 gcd(a,b,d)=\square\\-\frac{1}{2}d < b\leq \frac{1}{2}d
\\ X, Y\in \Z/n\Z
\end{array}}} [\sm a & b\\0&d\esm,(X,Y)](I-T)$$
$$+\sum_{\tiny{\begin{array}{cc}
ad=n^{2}, gcd(a,d)=\square\\ a,d>0, d\ \mbox{even}\\X,Y\in\Z/n\Z
\end{array}}}
[(\sm 0& -d\\ a& -\frac{1}{2}d\esm -\sm d& 0 \\ -\frac{1}{2}d & a
\esm ),(X,Y)].$$

The first sum on the right is $\equiv \hat{\mathcal{T}}_n^{\infty}(I-T),$ while
the second equals to
$$
\sum_{\tiny{\begin{array}{cc}
\alpha\beta=\frac{n^{2}}{2}, gcd(\alpha,\beta)=\square\\
\alpha,\beta>0\\
X,Y\in \mathbb{Z}/n\mathbb{Z}\end{array}}} [(\sm 0 & 2\beta\\ -\alpha &
\beta \esm -\sm 2\alpha & 0 \\ -\alpha & \beta\esm),(X,Y)]
$$
$$=\sum_{\tiny{\begin{array}{cc}
\alpha\beta=\frac{n^{2}}{2},
gcd(\alpha,\beta)=\square\\\alpha,\beta >0\\X,Y\in \mathbb{Z}/n\mathbb{Z}\end{array}}}
(S^2-I)[\sm 2\alpha & 0 \\ -\alpha & \beta\esm, (X,Y)] \equiv 0.$$
{\qed}
\subsection{Proof of Theorem \ref{Hecke1}-(1)}

To see $f|_{k,m}\mathcal{V}_{n}^{\infty}
\in J_{k,m}^{\int}(\Gamma(1))$ if  $f\in J_{k,m}^{\int}(\Gamma(1)) $
note that
$$\bigl((f|_{k,m}\mathcal{V}_{n}^{\infty})
\big|_{k,mn}\leb \SMA \alpha\beta\gamma\delta,
(\lambda,\mu)\rib\bigr)(\tau,z)$$
$$= n^{\frac{k}{2}-1}\biggl\{\sum_{ad=n\atop
b\MOD d}\biggl(f|_{k,m}\biggl[\MAT {\frac{1}{\sqrt{n}}}00{\frac{1}
{\sqrt{n}}},(0,0)\biggr]\biggr)|_{k,m}\leb \SMA{a/\sqrt{n}}{b/
\sqrt{n}}{c/\sqrt{n}}{d/\sqrt{n}},
(0,0)\rib\leb\SMA \alpha\beta\gamma\delta,
(\sqrt{n}\lambda,\sqrt{n}\mu)\rib\biggr\},$$ 
$$\mbox{ and } \, \bigl((f|_{k,m}\leb \SMA \alpha\beta\gamma\delta,
(\lambda,\mu)\rib)\big|_{k,m}\mathcal{V}_{n}^{\infty}\bigr)(\tau,z)$$
$$= n^{\frac{k}{2}-1}\biggl\{\sum_{ad=n\atop b\MOD d}
\biggl(f|_{k,m}\biggl[\MAT {\frac{1}{\sqrt{n}}}00{\frac{1}
{\sqrt{n}}},(0,0)\biggr]\biggr)|_{k,m}\leb\SMA \alpha\beta\gamma\delta, (\lambda,\mu)\rib\leb \SMA{a/\sqrt{n}}
{b/\sqrt{n}}{c/\sqrt{n}}{d/\sqrt{n}}, (0,0)\rib\biggr\}.
$$
So applying  the following equalities
$$\sum_{ad=l\atop b\MOD d}\leb \SMA{a/\sqrt{l}}{b/\sqrt{l}}{0}{d/\sqrt{l}},
(0,0)\rib\leb \SMA 1101, (0,0)\rib=\sum_{ad=l\atop b\MOD d}\leb
\SMA 1101, (0,0)\rib\leb \SMA{a/\sqrt{l}}
{b/\sqrt{l}}{0}{d/\sqrt{l}}, (0,0)\rib,$$
$$\sum_{ad=l\atop b\MOD d}\leb \SMA{a/\sqrt{l}}{b/\sqrt{l}}{0}{d/\sqrt{l}}, (0,0)
\rib\leb \SMA 1001, (0,\sqrt{l})\rib=\sum_{ad=l\atop b\MOD d}\leb \SMA 1001,
(0,a)\rib\leb \SMA{a/\sqrt{l}}{b/\sqrt{l}}{0}{d/\sqrt{l}}, (0,0)\rib$$
to $\biggl(f|_{k,m}\biggl[\MAT {\frac{1}{\sqrt{n}}}00{\frac{1}{\sqrt{n}}},(0,0)\biggr]\biggr)(\tau,z)$ we conclude our claim.

\subsection{Proof of Theorem \ref{Hecke1}-(2)}  By following similar way to the proof of Theorem \ref{Hecke2}-(2)
 we  find an explicit formula for $\tilde{\mathcal{V}}_n$ and the detailed proof will be skipped.

%%%%%%%%%%%%%%%%%%%%%%%%%%%%%%%%%%%%%%%%%%%

\end{document}